\documentclass[final]{siamltex}
\usepackage{graphics,graphicx,verbatim}
\usepackage{upquote}
\newlength\myverbindent 
\makeatletter
\def\verbatim@processline{%
 \hspace{\myverbindent}\the\verbatim@line\par}
\makeatother
\def\C{{\bf C}}
\def\Cp{{\bf C}^+}

\title{An algorithm for real and complex rational minimax approximation}

\author{Yuji Nakatsukasa\thanks{\texttt{nakatsukasa@maths.ox.ac.uk},
Mathematical Institute, University of Oxford, Oxford, OX2 6GG, UK.}\and
Lloyd
N.~Trefethen\thanks{\texttt{trefethen@maths.ox.ac.uk}, Mathematical
Institute, University of Oxford, Oxford, OX2 6GG, UK.}}

\begin{document}

\maketitle

\begin{abstract}
Rational minimax approximation of real functions on real intervals
is an established topic, but when it comes to complex functions
or domains, there appear to be no algorithms currently in use.
Such a method is introduced here, the {\em AAA-Lawson algorithm,}
available in Chebfun.  The new algorithm solves a wide
range of problems on arbitrary domains in a fraction of a second of
laptop time by a procedure consisting of two steps.  First, the
standard AAA algorithm is run to obtain a near-best approximation
and a set of support points for a barycentric representation of
the rational approximant.  Then a ``Lawson phase'' of iteratively
reweighted least-squares adjustment of the barycentric coefficients
is carried out to improve the approximation to minimax.

\end{abstract}


\pagestyle{myheadings}
\thispagestyle{plain}
\markboth{\sc Nakatsukasa and Trefethen}
{\sc Real and complex rational minimax approximation}

\section{Introduction}

Rational minimax approximation---the optimal approximation of a
function $f$ by a rational function $r$ of given degree on a given
domain in the supremum norm---is an old idea.  For real approximation
on a real interval, best approximations exist and are unique and
are characterized by an equioscillation condition.  Algorithms
appeared beginning in the 1960s~\cite{maehly,werner1,werner2},
and the problem became important for applications in the 1970s
with the development of digital signal processing~\cite{opp}.
A powerful implementation is available in the {\tt minimax} command
in Chebfun~\cite{chebfun,minimax}.  For complex functions or domains,
however, the situation is very different.  The theory developed by
Walsh in the 1930s shows that existence and especially uniqueness may
fail~\cite{walsh31,walsh}, and as for algorithms, there is not much
available apart from a pair of methods introduced by Ellacott and Williams (EW) (1976)
and Istace and Thiran (1993) based on earlier
work by Osborne and Watson~\cite{ow}, which, so far
as we are aware, are not in use today~\cite{ew,istace}.
(We have written our own EW code for comparisons; see later in this
section and also section~\ref{sec-complex}.)
This is a striking gap, since rational approximations are of
growing importance in systems theory and
model order reduction~\cite{antoulas,loewner,bcow,chdo},
electronic structure calculation~\cite{pexsi,moussa}, low-rank
data compression~\cite{becktown,liwhite}, computational complex
analysis~\cite{gt1}, and solution of partial differential
equations~\cite{gt2,gt3}.

The aim of this paper is to introduce a new algorithm for complex
rational minimax approximation together with a software implementation.
Our ``AAA-Lawson'' algorithm combines the rational barycentric AAA
algorithm of~\cite{aaa} with an iteratively reweighted least-squares
(IRLS) iteration inspired by Lawson's algorithm~\cite{lawson} but
in a nonlinear barycentric context.  It works on discrete domains,
typically containing hundreds or thousands of points to approximate
a continuum, which may take all kinds of forms including Jordan
regions, unions of Jordan regions, regions with holes, intervals,
unbounded domains, clouds of random points, and more.  Being based
on a barycentric rational representation with greedy selection of
support points, it inherits the exceptional numerical stability of
AAA and is able to handle even very difficult cases with exponentially
clustered poles.  Experiments show that for a wide range of problems,
the method converges in a fraction of a second on a laptop to an
approximation with an error within a few percent of the minimax value.

The version of the AAA-Lawson algorithm described here was introduced
in the {\tt aaa} command of Chebfun in August, 2019, and we hope its
easy availability may open up a new era of exploration of complex
rational minimax approximation.  For example, Ellacott and Williams
list minimax errors for 29 different approximation problems in Tables~1
and~2 of their paper, each given to 3 digits of accuracy~\cite{ew}.
With Chebfun {\tt aaa}, all of these results can be computed by a few
lines of code in a total time of less than 2 seconds on a laptop.
(Twelve of the Ellacott--Williams numbers turn out to be correct
to all three digits, with the rest having small anomalies mainly
associated with discretization of a continuum by too few points.
As we shall discuss in section~\ref{sec-complex}, however, the
EW method works for only a limited range of problems.)

For complex polynomial approximation, more computational
possibilities are available than in the rational case,
including~\cite{bdm,fm,gr,ko,opfer,tang}.  The ``complex Remez
algorithm'' of Tang is particularly appealing~\cite{tang}.  Similarly
there are a number of non-minimax complex rational approximation
algorithms, including vector fitting~\cite{gs}, RKFIT~\cite{bg}, the
Loewner framework~\cite{loewner}, IRKA~\cite{gab}, and AGH~\cite{agh},
as well as the AAA algorithm that is our own starting point~\cite{aaa}.
Most of these methods apply to vector or matrix as well as scalar
approximation problems, whereas the AAA-Lawson method has only been
developed so far for scalars.  To extend it, one could perhaps adapt
some of the methods proposed for AAA in~\cite{lietaert}.

The possibility of a AAA-Lawson algorithm was first mentioned in
the original AAA paper~\cite{aaa}, and it was developed further for
part of the initialization process for the Chebfun {\tt minimax}
command~\cite{minimax}.  However, in these projects the power of
AAA-Lawson for general minimax approximation did not become fully
apparent, for a number of reasons.  One was that AAA approximations
are usually computed with the degree not specified but adaptively
chosen to get down to nearly machine precision, and in this setting,
AAA-Lawson will usually fail (it is trying to improve a result that
is already near the limit of precision).  Another is that much of
our attention was on real intervals, where both AAA and AAA-Lawson
are least robust.  A third was that we did not fully appreciate the
crucial importance of choosing approximation grids exponentially
clustered near corners and other singular points,
where poles of rational approximations will
be exponentially clustered.  Finally, in those experiments we were
not including the support points themselves in the matrix associated
with the IRLS problem (see eq.~(\ref{special})), an omission that
led to failure in some cases.

\begin{figure}
\begin{center}
\includegraphics[scale=.85]{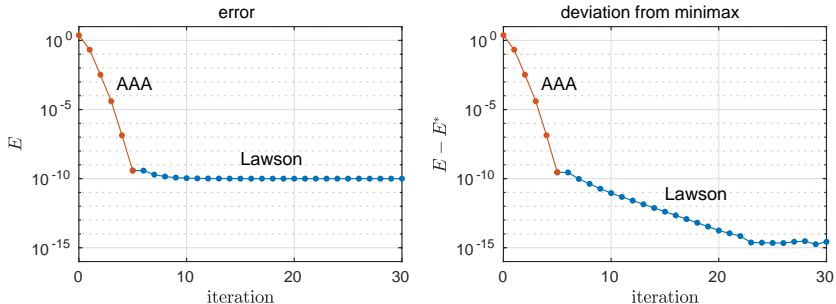}
\end{center}
\caption{\label{fig1}The two phases of the AAA-Lawson algorithm, illustrated
here for degree $5$ approximation
of $e^z$ on the unit disk.  The AAA phase achieves rapid convergence
to a near-minimax approximant.  This is then improved to minimax by
a linearly convergent Lawson iteration.}
\end{figure}

We close this introduction with Fig.~\ref{fig1},
illustrating behavior
of the algorithm in a typical problem (the first example
of section~\ref{sec-complex}).
The first, AAA phase rapidly finds a near-minimax approximation,
and this is improved to minimax in the second, Lawson phase.

\section{\label{theory}Existence, uniqueness, characterization, and convergence}
A rational function is a function of a complex variable that can
be written in the form $r(z) = p(z)/q(z),$ where $p$ and $q$ are
polynomials.  We say that $r$ is of {\em type $(m,n)$} for some $m,
n \ge 0$ if it can be represented with $p$ of degree at most $m$
and $q$ of degree at most $n$.  If $m=n$, the setting of this paper,
we also say that $r$ is of {\em degree $n$}, and we denote by $R_n$
the set of rational functions of degree $n$.  A rational function is a
map from $\Cp$ to $\Cp$, where $\Cp$ is the extended complex plane $\C
\cup\{\infty\}$, and if $r\in R_n$ is not a constant, then it takes
each value in $\Cp$ at most $n$ times, counted with multiplicity.
It is the function that is the fundamental
object, not any particular representation
of it, and if a representation has isolated points corresponding to
quotients $0/0$ or $\infty/\infty$, we define the values there by
limits from neighboring points.

Let $Z\subseteq \C$ be nonempty, let $\|\cdot\|$ be the supremum norm
on $Z$, and let $f$ be a complex continuous function (not necessarily
analytic) defined on $Z$.  Our approximation problem is to find
rational functions $r$ such that $\|f-r\|$ is small.  If $r^*\in R_n$
satisfies $\|f-r^*\| = E^* := \inf_{r\in R_n} \|f-r\|$, then $r$ is
a {\em best\/} or {\em minimax\/} approximation to $f$ of degree $n$.
Even if there is no minimax approximation, we will speak of the {\em
minimax error} $E^*$, which may be any number in the range $[\kern
.3pt 0,\infty]$.

Polynomial best approximations of a given degree $n$ always exist,
and if $Z$ is a compact set with at least $n+1$ points, they are
unique~\cite{tonelli}.  They can be characterized by a condition due
to Kolmogorov~\cite{kol,opfer,singer,thiran}, and if $Z$ is a closed real
interval of positive length and $f$ is real on $Z$, there is a simpler
and more famous characterization by equioscillation of the error
$(\kern .5pt p-f)(x)$ between $\ge n+2$ alternating extreme points.
In this case $E^*$ decreases exponentially as $n\to\infty$ if and
only if $f$ is analytic on $Z$~\cite[chapter 8]{atap}, and the same
result generalizes to any compact set $Z\subseteq \C$~\cite{walsh}.

The theory of complex rational best approximation, which begins with
a 1931 paper by Walsh~\cite{walsh31,walsh}, is not so satisfactory.
First of all, best approximations need not exist.  For example,
there is no degree $1$ best approximation to data values $a,b,b$
with $a\ne b$ on any set $Z$ with three distinct points, for in
such a case we have $E^* = 0$ but $\|f-r\|>0$ for any choice of
$r$, since a nonconstant function $r\in R_1$ can only take each
value once.  However, Theorem~3 of~\cite{walsh31} asserts that
existence is assured if $E^*<\infty$ and $Z$ has no isolated points.

Concerning uniqueness, there is one main positive result:
if $Z$ is a closed real interval of positive length and $f$
is real, then a best {\em real\/} rational approximation in
$R_n$ exists and is unique and characterized by an error curve
that equioscillates between sufficiently many extreme points.
Without the restriction that $r$ is real, however, uniqueness is
not assured~\cite{lungu,sv,williams}.  For example, there
are complex approximations to $|x|$ on $[-1,1]$ whose error is
less than the value $1/2$ achieved by the best real approximation,
from which it follows by the symmetry of complex conjugation that
the complex best approximation cannot be unique (see Exercise
24.3 of~\cite{atap}).  Examples of nonuniqueness have also been
investigated on the unit disk~\cite{gt}.

When it comes to characterization of rational best approximants,
the Kolmogorov condition diminishes to a necessary condition for
local optimality: for a candidate approximation to be locally
optimal, it must be a stationary point with respect to certain
local linear perturbations.  Discussions can be found in a number
of sources, including~\cite{gutk1,ruttan,singer,williams}, and we recommend in
particular the papers~\cite{istace,thiran} by Istace and Thiran.
Sufficient conditions, and conditions for global optimality,
are mostly not available, though there are some results in~\cite{ms}
and~\cite{ruttan}.  We shall say more on these subjects in
sections~\ref{sec-complex} and~\ref{steps}.

These observations are daunting.  However, although it is a
fascinating mathematical challenge to elucidate the properties
of {\em best} approximations, what matters for most applications
is that we are able to compute {\em good\/} ones.  An example is
provided by the ``lightning Laplace solver'' paper \cite{gt3}, which
presents far-reaching theorems about root-exponential convergence
of rational approximations for solutions to Laplace problems with
boundary singularities.  The approximations are not minimax, but
still they lead to a very fast Laplace solver.

A major focus of the theoretical literature of rational approximation
is the problem of {\em approximability,} the determination
of necessary and sufficient conditions to ensure that $E^*\to 0$
as $n\to \infty$ in approximation of a function $f$ on a set
$Z\subseteq \C$.  If $Z$ is compact with at most a finite number
of holes and $f$ is analytic on $Z$, then exponential decrease of
$E^*$ to $0$ was established by Runge in 1885~\cite{runge},
but what if $f$ is merely analytic in the interior?  Here we cannot
expect exponential convergence, but according to {\em Vitushkin's
theorem}~\cite{gaier,gamelin,zalcman}, the generalization to
rational approximation of Mergelyan's theorem for polynomials, we
still get $E^*\to 0$.  And what if there are infinitely many holes?
Vitushkin's theorem gives technical conditions for this case too.
But such questions are a long way from most applications of rational
approximation, where the whole point is to exploit circumstances
in which $E^*\to 0$ very fast.

\section{\label{sec-alg}The AAA-Lawson algorithm}

Let $n+1$ distinct {\em support points\/} $t_0,\dots, t_n\in\C$ be fixed
for some $n\ge 0$, and let
$\ell$ be the {\em node polynomial}
\begin{equation}
\ell(z) = \prod_{k=0}^n (z-t_k),
\label{node}
\end{equation}
which is monic and of degree $n+1$.  If
$\alpha_k, \beta_k\in\C$ are arbitrary complex numbers, $0\le k \le n$,
with at least one $\beta_k$ being nonzero,
then the quotient of partial fractions
\begin{equation}
r(z) = {n(z)\over d(z)} = \sum_{k=0}^n {\alpha_k\over z-t_k} \left/
\sum_{k=0}^n {\beta_k\over z-t_k} \right.
\label{bary}
\end{equation}
is obviously a rational function of degree $2n+1$, since
the numerator $n$ and denominator $d$ are rational functions of type $(n,n+1)$ and
the denominator is not identically zero.
However, by multiplying both $n$ and $d$ by $\ell$,
we see more sharply that $r$ is of degree $n$.  The expression
(\ref{bary}) is a {\em barycentric representation} for $r$~\cite{bt}.

Conversely, regardless of the choice of the support points, every
degree $n$ rational function $r$ can be written in the form 
(\ref{bary}).
The following theorem and the first proof are adapted from~\cite{aaa}.

\medskip

\begin{theorem}[Rational barycentric representations]
Let $t_0,\dots, t_n$ be an arbitrary set of distinct complex
numbers.  As $\alpha_0,\dots ,\alpha_n$ and
$\beta_0,\dots ,\beta_n$
range over all complex values, with at least one
$\beta_k$ being nonzero, the functions $(\ref{bary})$
range over the set of all
rational functions of degree $n$.
\end{theorem}

\smallskip

\begin{proof}
As just observed, any quotient
(\ref{bary}) is a rational function $r$ of degree $n$.
Conversely, suppose $r$ is a rational function of degree $n$,
and write $r = p/q$ where $p$ and $q$ are
polynomials of degree at most $n$.
Then $q/\ell$ is a rational function with a zero at $\infty$
and a simple pole at each point $t_k$, or no pole
at all if $q(t_k)=0$.  Therefore $q/\ell$
can be written in the partial fraction form $d$
as in (\ref{bary}) (see p.~553 of~\cite{henrici}).
Similarly $p/\ell$ can be written in the form $n$.
\end{proof}

{\em Alternative proof.}
Writing $r=p/q$ again, we note that it is enough to show
that coefficients $\{\alpha_k\}$ and $\{\beta_k\}$ exist such 
that $p = n\kern .3pt \ell$ and $q = d \kern .3pt \ell$ in
(\ref{bary}).
Now $d\kern .3pt \ell$ is a linear combination
with coefficients $\beta_0,\dots, \beta_n$
of $n+1$ monic polynomials of degree $n$, which 
are linearly independent since they have different
sets of roots.  Thus $q$ can be written (uniquely)
as $d\kern .3pt \ell$, and similarly for $p = n\kern .3pt \ell$.
\endproof

\medskip

The second of these proofs shows that there is
a one-to-one correspondence between sets of coefficients $\{\alpha_k\}$
in a barycentric representation (\ref{bary}) and polynomials
$p$ in a quotient representation $p/q$, and likewise for
$\{\beta_k\}$ and $q$.  Thus we see that the
barycentric representation is unique
to exactly the same degree as the quotient representation $p/q$: unique
up to a multiplicative constant if $r$ has degree $n$ but
not $n-1$, with further nonuniqueness if $r$ is of degree
$n-1$ or less.

Rational barycentric formulas
with independent coefficients $\alpha_k$ and $\beta_k$ are not well known.
Traditionally, barycentric formulas are used in ``interpolatory mode,''
where function values $\{f_k\}$ are given and weights are chosen corresponding
to $\alpha_k/\beta_k = f_k$ (and $\beta_k \ne 0$), yielding
$r(t_k) = f_k$ for each $k$~\cite{bt,fh,aaa,sw}.  To
work with arbitrary rational functions, however, with a complete decoupling of
support points from approximation properties, one needs the
``noninterpolatory'' or ``alpha-beta'' mode (\ref{bary}).
Ultimately the $\alpha_k$ and $\beta_k$ are devoted
to approximation and the $t_k$ to numerical stability.

The AAA-Lawson algorithm consists of two steps.  We assume that a
discrete domain $Z$ and a set of corresponding function values $F=f(Z)$
have been given, together with a degree $n$.

\medskip

{\em {\rm (I)} Run the AAA algorithm to get a rational approximant
$r_0\approx f$ of degree $n$ and a set of support points
$t_0,\dots, t_n$.}

\smallskip

{\em {\rm (II)} Carry out a linearized barycentric
Lawson iteration until a termination condition is reached.}

\medskip

\noindent Step (I) utilizes (\ref{bary}) in interpolatory
mode, with the support points chosen one after another in a
greedy manner.  Typically $\|f-r_0\|$ is within about an order of
magnitude of the minimax error, but
since $r_0$ interpolates the data
at $n+1$ points, one cannot expect it to be the optimal approximant.
The details of (I) are presented in~\cite{aaa}, and we shall
not repeat them here.  What remains is to describe step (II), which
switches to noninterpolatory mode.

Let $Z = (z_j), $ $1\le j \le M$ be the sample set, interpreted
as a column vector, and let
$F = (f_j), $ $1\le j \le M$ be the corresponding vector of
function values to be matched.  Let $\alpha$ and $\beta$ be the
coefficient vectors
$(\alpha_0,\dots,\alpha_n)^T$ and
$(\beta_0,\dots,\beta_n)^T$, with $\gamma$ defined
as their concatenation $\gamma = [\kern .7pt \alpha ;\,  \beta\kern .7pt ]$.
Our aim is to solve the minimax problem
\begin{equation}
\min_\gamma \,
\max_j \left| f_j - \sum_{k=0}^n {\alpha_k\over z_j-t_k} \left/
\sum_{k=0}^n {\beta_k\over z_j-t_k} \right. \right|.
\label{bary2}
\end{equation}
The barycentric Lawson
idea is to achieve this by solving a sequence of
iteratively reweighted least squares (IRLS) problems 
based on the linearization of (\ref{bary2}),
\begin{equation}
\min_{\gamma,\, \|\gamma\|_2=1} \;
\sum_{j=1}^M {}' \; w_j \left(f_j \sum_{k=0}^n {\beta_k\over z_j-t_k} -
\sum_{k=0}^n {\alpha_k\over z_j-t_k}\right)^2 ,
\label{barylin}
\end{equation}
where at each step, $W = (w_j)$, $1\le j \le M,$ is a vector
of weights $w_j\ge 0$.
Note the prime symbol on the summation sign.  This signifies that
special treatment is applied at the $n+1$ sample points $z_j$
that coincide with a support point $t_k$ for some $k = k_j$.  At
these points the quantity in parentheses in (\ref{barylin})
would be infinite, and instead, in
the spirit of L'H\^opital's rule, these terms of
the sum are replaced by 
\begin{equation}
 w_j \Bigl(f_j \beta_{k_j}^{} - \alpha_{k_j}^{}\Bigr)^2 .
\label{special}
\end{equation}

Equation (\ref{barylin}) is a routine problem of numerical
linear algebra, which can be written in matrix form as
\begin{equation}
\min_{\gamma, \,\|\gamma\|_2=1} \;
\bigl\| \kern 1pt \hbox{diag}\kern .3pt (W^{1/2})
\kern 1pt \bigl[\kern 1pt C,
-\hbox{diag}\kern .3pt (F) \kern .7pt C\kern 1pt
\bigr ]\kern 1.2pt  \gamma\kern 1pt  \bigr\|_2,
\label{matrixform}
\end{equation}
where $C$ is the {\em Cauchy matrix\/} with entries
$c_{ij} = 1/(z_j-t_k)$ except in the $n+1$ special rows.
This is a minimal singular value problem
involving a matrix of size $M\times (2n+2)$,
and it can be solved on a laptop in on the order of a millisecond
for typical values of, say, $M=1000$ and $n=20$.

From one IRLS step to the next, $W$ is updated by the formula
\begin{equation}
w_j^{\hbox{\scriptsize(new)}} = w_j |e_j|,
\label{update}
\end{equation}
where $e_j$ is the quantity inside absolute values in (\ref{bary2}), 
i.e., the current nonlinear error at $z_j$.  (For the $n+1$
special values of $j$, $e_j = f_j - \alpha_{k_j}^{}/\beta_{k_j}^{}$.)
For convenience,
and floating-point arithmetic, we then renormalize the
weights at each step so that their maximum is~$1$.

The IRLS idea originated with Lawson in 1961~\cite{lawson} for
linear minimax approximation, and has subsequently been analyzed and
generalized by
a number of authors beginning with Cline, Rice, and Usow~\cite{cline,rice,ru}.
Rice proved convergence at a linear rate
for real approximation under natural assumptions~\cite{rice},
and Ellacott and Williams pointed out that the same proof extends to
complex approximation~\cite{ew}.  IRLS algorithms have also taken
on importance for other kinds of linear $L^p$ approximation, particularly
the case $p=1$ of interest in data science~\cite{ddfg,osborne,watson}.
However, apart from a (non-barycentric) attempt with
limited success in~\cite{cooper}, AAA-Lawson is the first IRLS algorithm we know of for
nonlinear approximation.

This completes our description of the core idea of the AAA-Lawson
algorithm, but three questions remain.
(i) How do we terminate the iteration? (ii) What can be proved about
convergence?  (iii) What steps can be taken to make convergence
faster or more reliable in troublesome cases?  Even for linear
approximation, these are nontrivial matters, and the
nonlinear case brings additional difficulties.  Chebfun's
answer to (i) is simple: by default it takes 20 Lawson steps.  We shall
discuss (ii) and (iii) in section~\ref{steps}.

\section{\label{sec-complex}Numerical examples, complex}

In this section we present fourteen examples of complex minimax
approximations, grouped into pairs for convenience.  Each example
is represented by three images in the complex plane, the first
showing the domain $Z$ and the second and third showing the error
$r(Z)-f(Z)$ in AAA and AAA-Lawson approximation.  All computations
were done in Chebfun in the default mode, and laptop timings
are printed at the tops of the figures.  The codes of this section
and the next are available in the supplementary materials.

For comparison, we have also solved all the problems of this section by
the Ellacott--Williams (EW) method from~\cite{ew} (our own implementation),
getting correct results in about half of the cases.
We find that when EW is successful, it is typically 5--100 times slower
than AAA-Lawson since each step requires the iterative
solution of a linear optimization problem.
(Its asymptotic convergence rate is usually better, however, so the
timings get closer if
minimax errors are required to many digits of accuracy.  Also, the
Istace--Thiran algorithm~\cite{istace}, which we have not implemented, is likely 
to be faster than EW.)
The problems where EW is successful are those
involving domains not too far from the unit circle and without
singularities on the boundary, as in
Figs.~\ref{expzpair}, \ref{squarepair}, and~\ref{arcpair} below.
For other problems, such as those
of Figs.~\ref{scpair} and \ref{pair7} and the NICONET
problem of Fig.~\ref{lastpair}, it generally 
fails to find the minimax solution, which we attribute to
its reliance on $p/q$ rather than barycentric representation of
rational functions.
Explanations of the great difference in stability between
these two representations can be found
in~\cite[sec.~2]{minimax}
and~\cite[sec.~11]{aaa}.

\begin{figure}
\begin{center}
\includegraphics{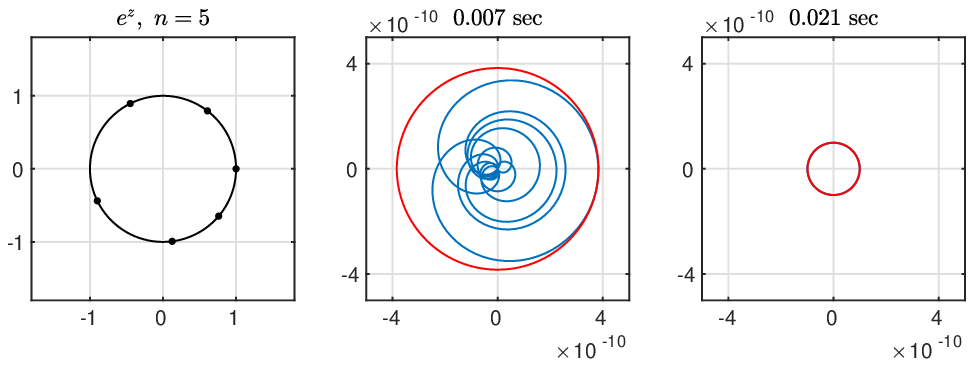}
\smallskip
\includegraphics{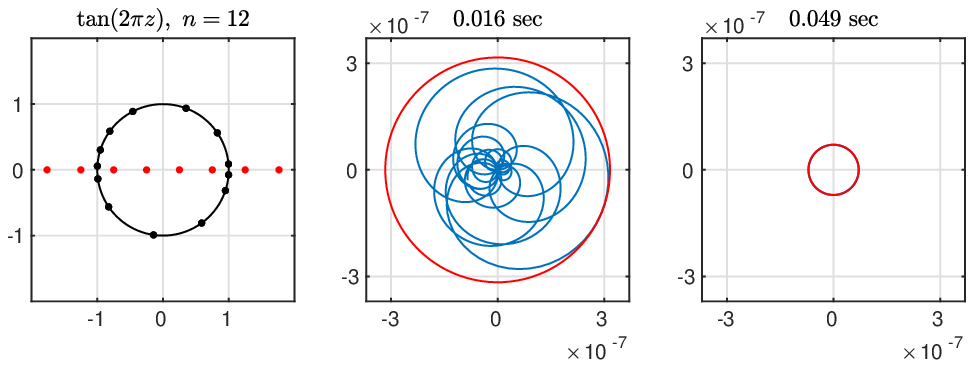}
\end{center}
\caption{\label{expzpair}Approximation on the unit circle of
the analytic function $e^z$ and the meromorphic function
$\tan(2\pi z)$.  The middle and right plots show error curves
for AAA and AAA-Lawson approximation, respectively, with red
circles marking the maximum errors.  The
minimax error curves are nearly circular (invisible
here since they lie under the red circles), with winding
numbers $2n+1 = 11$ and $2n+1-8 = 17$.
Red dots mark poles of the AAA-Lawson approximation.}
\end{figure}

Figure~\ref{expzpair} begins with the basic example of $e^z$ on
the unit circle, discretized by 500 equispaced points.  With its
usual great speed, AAA finds a near-best approximation for $n=5$
with error 3.83e-10.  The black dots on the circle mark the six
support points the algorithm selects.  Continuing with the same support
points but now in noninterpolatory ``alpha-beta'' mode, AAA-Lawson
improves the approximation to close to
minimax, with error 9.944364e-11.
By the maximum modulus principle for analytic functions, these
maximal errors on the circle are also the maximal errors over the
whole disk.  Note that the error curve appears to be a perfect circle
(of winding number $2n+1 = 11$, though this cannot be distinguished
in the figure).  This near-circularity effect was first identified
in~\cite{ugrad} and then investigated for polynomial approximation
in~\cite{nearcirc} and rational approximation in~\cite{cf}.
The error curve cannot be exactly circular (this would imply that the
function being approximated was rational), but as shown in~\cite{cf},
it comes spectacularly close, varying in radius for this example, we
estimate via Theorem 6.3 of~\cite{cf}, by
less than one part in $10^{12}$.  This effect led to the theory of
Carath\'eodory--Fej\'er (CF) approximation~\cite{htg,cf,atap,vandeun},
which establishes the lower bound $E^*\ge \sigma_{n+1}$, where $\{\sigma_k\}$
are the singular
values of the infinite Hankel matrix of Taylor coefficients
$a_1, a_2, \dots = 1, 1/2\kern 1pt !, \dots.$  Here the relevant
value is $\sigma_6 = \hbox{9.944144081e-11}$.

The second example of Figure~\ref{expzpair} is $\tan(2\pi z)$
for $n=12$ in 1000 points of the unit circle.  This function is
meromorphic but not analytic in the unit disk.  Again we get a
nearly-circular error curve, whose winding number is not 25 but 17
because of the four poles in the disk.  Here AAA-Lawson improves the
error from 3.16e-7 to 7.08e-8.  The red dots in the left image mark
poles of the AAA-Lawson approximation.  The poles inside
the circle match the poles $\pm 1/4$ and $\pm 3/4$ of
$\tan(2\pi z)$ to 13 digits of accuracy.  We can explain this
by noting that these poles can be determined by
certain contour integrals of the boundary data~\cite[sec.~4]{akt},
and since $r$ matches
$\tan(2\pi z)$ to many digits on the boundary, the contour
integrals must match too.
The poles of $r$ outside the circle are at $\pm 1.250011$,
$\pm 1.7638$, $\pm 2.6420$ and $\pm 7.3844$.  (In the first row of
this figure, no red dots appear because the poles are off-scale.
Their positions in the case of Pad\'e approximations were
investigated by Saff and Varga~\cite{svpoles}.)

\begin{figure}
\begin{center}
\includegraphics{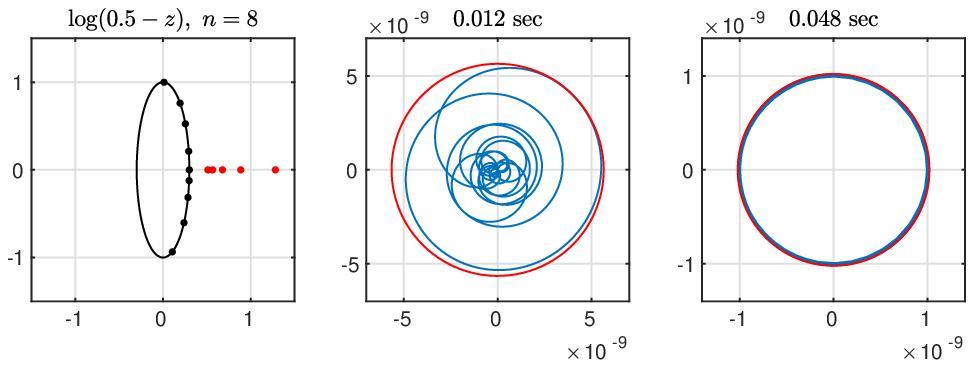}
\smallskip
\includegraphics{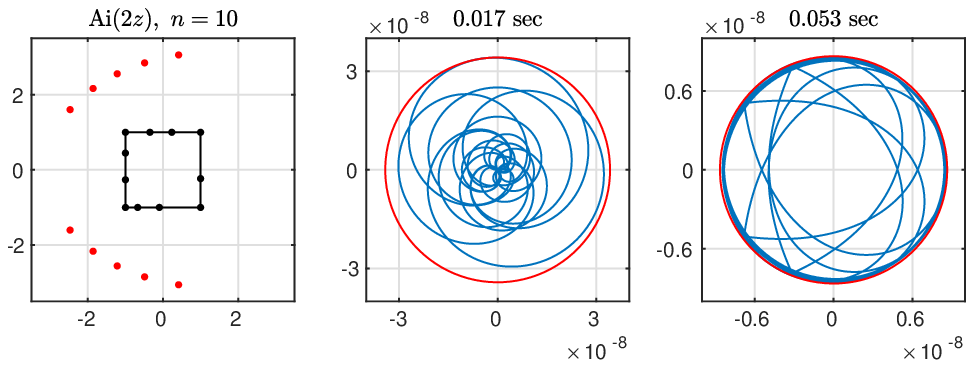}
\end{center}
\caption{\label{squarepair}Approximations on an ellipse and a square.
The near-circularity effect appears again, though on the square, the four
corners persist.
Here and in most of figures to follow, the axis scales are different
for the AAA and AAA-Lawson plots.}
\end{figure}

Figure~\ref{squarepair} shows approximations on two noncircular
domains.  In the first row, $\log(0.5-z)$ is approximated in 2000
points on an ellipse of half-height 1 and half-width $0.3$.  Note how
the poles of the approximation line up along the branch cut, a
phenomenon analyzed for Pad\'e approximations by Stahl~\cite{stahl}.
It is also interesting to see that all the support points chosen
by AAA lie on that side of the ellipse.  The second row shows
approximation of the Airy function $\hbox{Ai}(2z)$ in 4000 points
on the boundary of the unit square, 1000 points in a Chebyshev
distribution on each side.  The error curve, with winding number
$2n+1 = 21$, is nearly-circular along most of its length, while
retaining the four corners associated with the square.

\begin{figure}
\begin{center}
\includegraphics{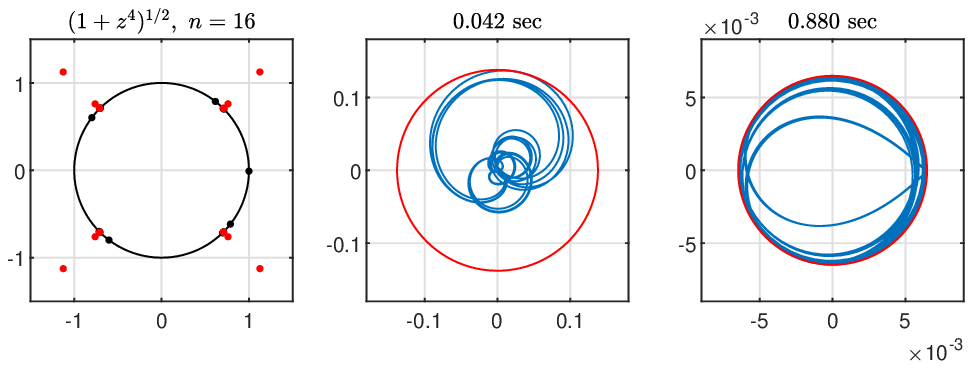}
\smallskip
\includegraphics{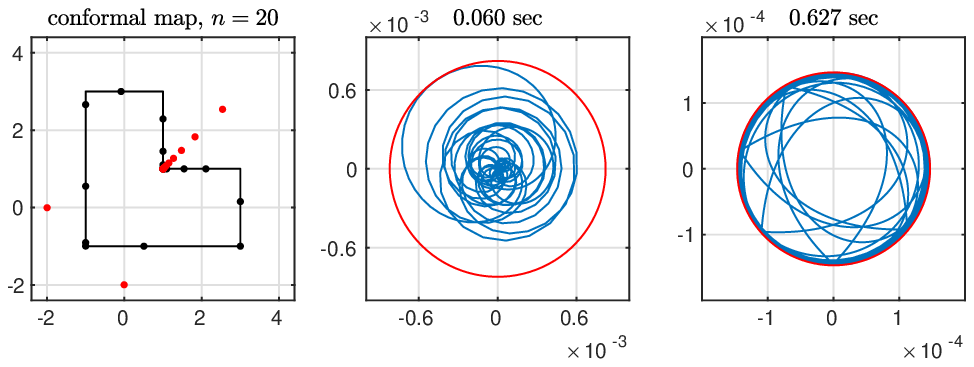}
\end{center}
\caption{\label{scpair}Two approximation problems with singularities on the
boundary.  The second row,
following~{\rm \cite{gt1}}, is the inverse of a Schwarz--Christoffel
conformal map.  Because of the prevalence of corner
singularities, rational approximations can be a powerful
tool in numerical conformal mapping.}
\end{figure}

Figure~\ref{scpair} turns to problems with singularities on the
boundary, where rational functions have their greatest power
relative to polynomials, achieving root-exponential convergence
as $n\to\infty$ by means of poles exponentially clustered near
the singularities~\cite{gt1,gt2,gt3,newman}.  In the first row,
$(1+z^4)^{1/2}$ is approximated to degree $n=16$.  The AAA-Lawson
approximation improves the error from 1.38e-1 to 6.49e-3, with poles
lying along branch cuts near each of the four singularities at radii
1.00046, 1.0085, 1.075, and 1.59.  For successful computation of
approximations with clustered poles like this, it is important that
the sample grid be clustered too, and in this case the sample points
on the unit circle were placed at angles $(\pi/4)\cdot\hbox{{\tt
tanh(linspace(-12,12,1000))}}$ together with their rotations by
angles $\pi/2$, $\pi$, and $3\pi/2$.  Note that there are four square
roots in this function, hence four right angles in the error curve,
but these appear as one because they lie on top of one another.

The second row of Fig.~\ref{scpair} shows degree 20 approximation
of an analytic function representing a conformal map of an L-shaped
region onto the unit disk, which has a $z^{2/3}$ type of singularity
at the reentrant corner.   Each of the six sides has sample points
with a distribution controlled by {\tt tanh(linspace(-12,12,1000))}.
In~\cite{gt1} it was shown that AAA rational approximations of
conformal maps of polygons can be 10--1000 times more efficient to
evaluate than the standard method of Driscoll's Schwarz--Christoffel
Toolbox~\cite{toolbox}.  From Figure~\ref{scpair} we see that
even better approximations are available with AAA-Lawson, which
improves the accuracy of the approximation in this case from 8.21e-4
to 1.57e-4.

\begin{figure}
\begin{center}
\includegraphics{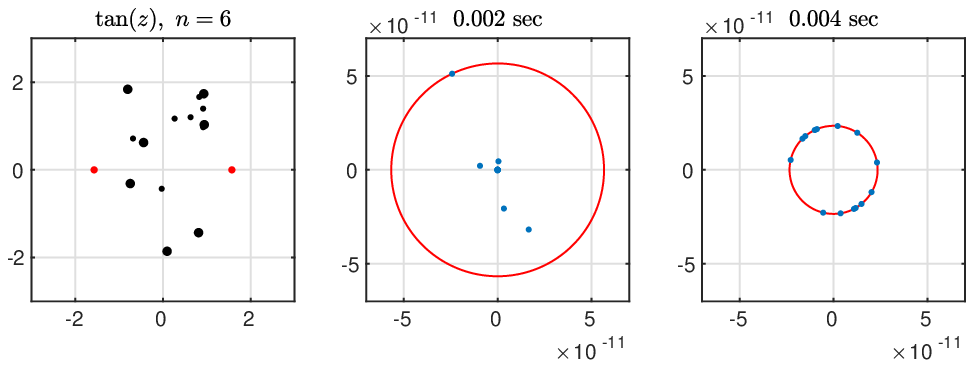}
\smallskip
\includegraphics{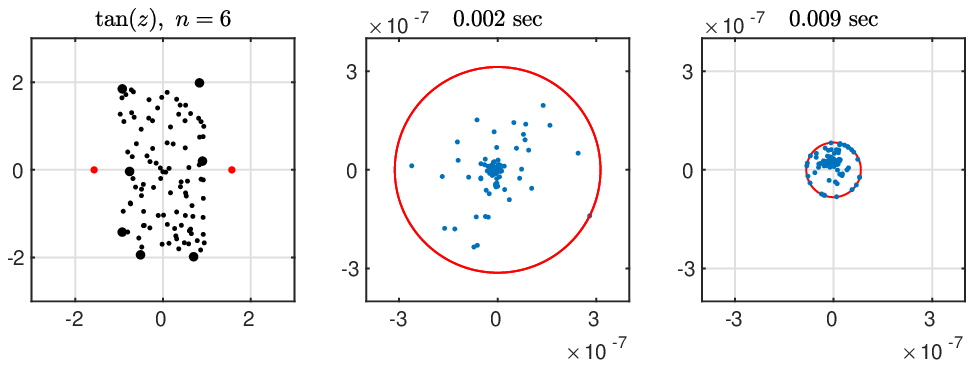}
\end{center}
\caption{\label{randpair}Approximation of\/ $\tan(z)$ 
at $14$ and $100$ random points in a rectangle in $\C$.  In the first case, with
just $2n+2$ sample points, the minimax error is attained at every one.}
\end{figure}

Figure~\ref{randpair} moves from essentially continuous domains
to discrete ones consisting of random points in a rectangle.
A rational function of degree $n$ could generically interpolate
$2n+1$ data values exactly.  Thus the first nontrivial fit occurs
with $2n+2$ data values, and this is shown in the first row of the
figure, with $n=6$ and 14 sample points.  As expected, the minimax
error is attained at all 14 points.  The second row increases
the number of sample points to 100, and now the maximum error, which
is 10,000 times larger, is
attained at 20 rather than 14 of them.  (This is not evident with
the calculation as run with the Chebfun default number of 20 Lawson
steps, but emerges if a few hundred Lawson steps are taken to give
convergence to more digits of accuracy.)

\begin{figure}
\begin{center}
\includegraphics{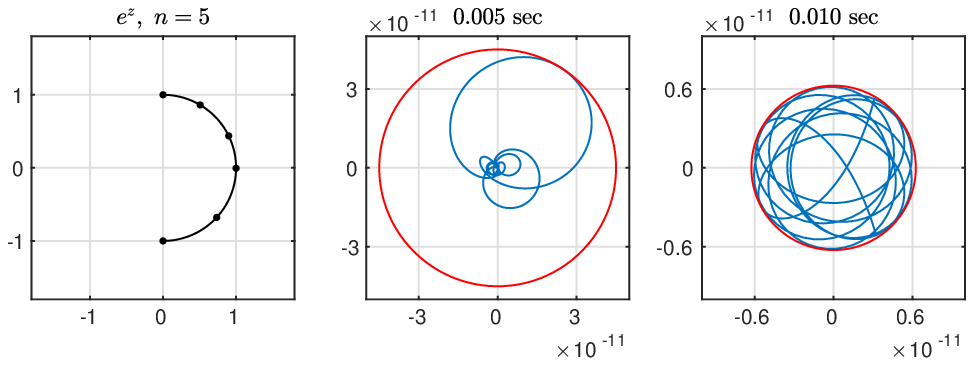}
\smallskip
\includegraphics{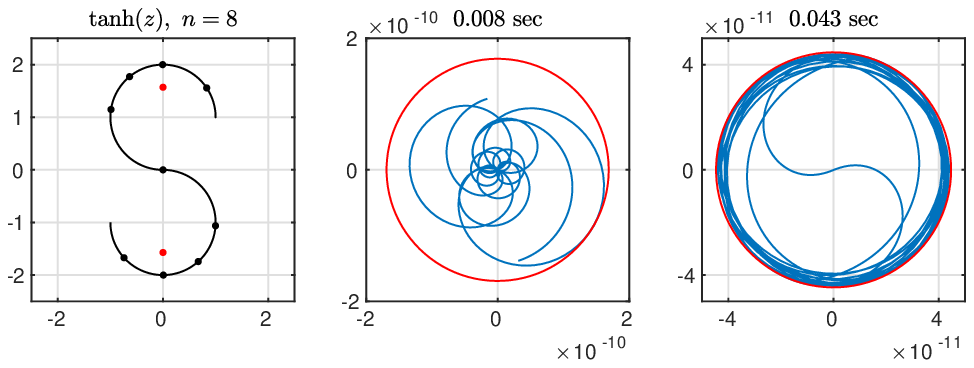}
\end{center}
\caption{\label{arcpair}Approximations on complex arcs.}
\end{figure}

Figure~\ref{arcpair} shows two approximations on domains that are
just arcs, a semicircle and an S-shape, both represented with 500
points in a Chebyshev distribution along each semicircular piece.
Figure~\ref{lastpair} shows, first, an approximation on the unit
circle of a function with an essential singularity in the disk,
and second, the clamped beam example from the NICONET model
order reduction collection~\cite{chdo}, which was also considered
in~\cite{aaa}.  Here the approximation domain is the imaginary axis,
which is discretized by 2000 points logarithmically spaced between
$0.01i$ and $100\kern .3pt i$ together with their complex conjugates.
The function to be approximated is defined via the resolvent of a
$348\times 348$ matrix whose eigenvalues are in the left half-plane,
making it analytic in the right half-plane.  Note that in this
example, AAA-Lawson achieves reduction of the error by a factor
of about 4, from 6.15 to 1.49.  (This is $0.03\%$ accuracy, for
the function being approximated takes values as large as 4550.)
In an application of model order reduction, such an improvement
might be significant~\cite{loewner,bcow}.

\begin{figure}
\begin{center}
\includegraphics{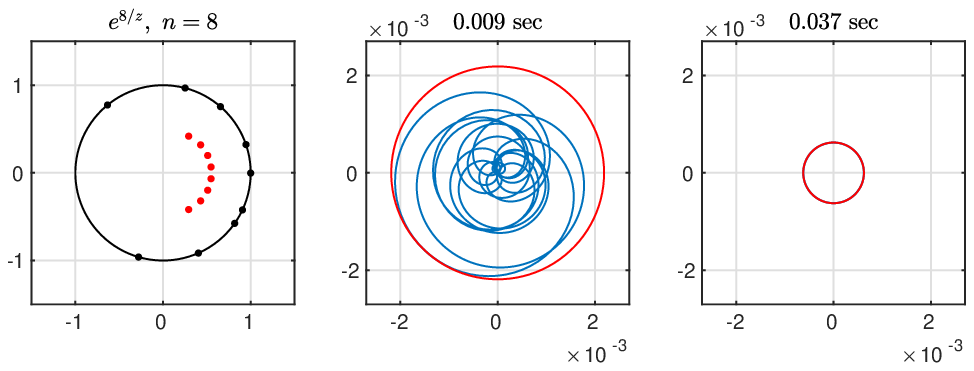}
\smallskip
\includegraphics{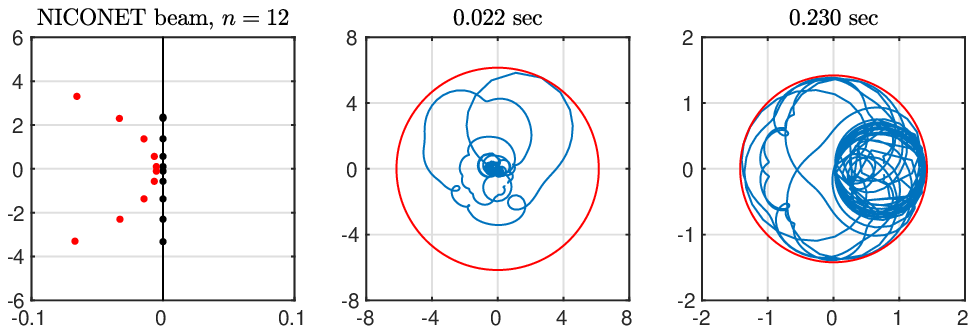}
\end{center}
\caption{\label{lastpair}The example of the first row has an essential singularity
in the unit disk;
all the poles fall inside the circle and the winding number is $-17$.
The second is the NICONET beam model order reduction
example of\/ {\rm \cite{chdo}}, defined on the
imaginary axis via the resolvent of a $348\times 348$ matrix.}
\end{figure}

\begin{figure}
\begin{center}
\includegraphics{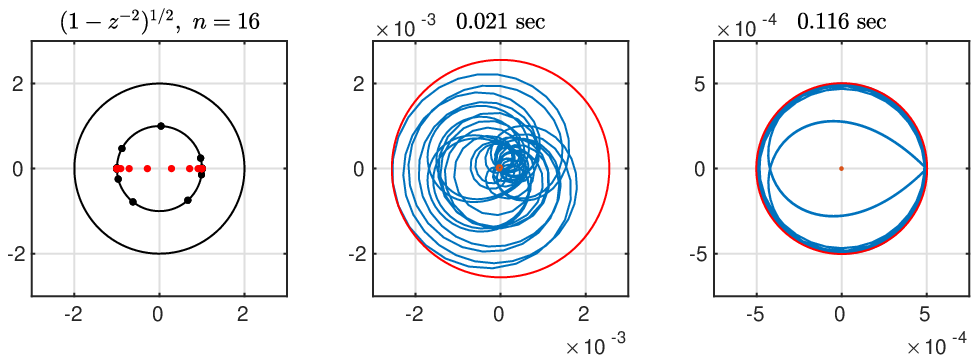}
\smallskip
\includegraphics{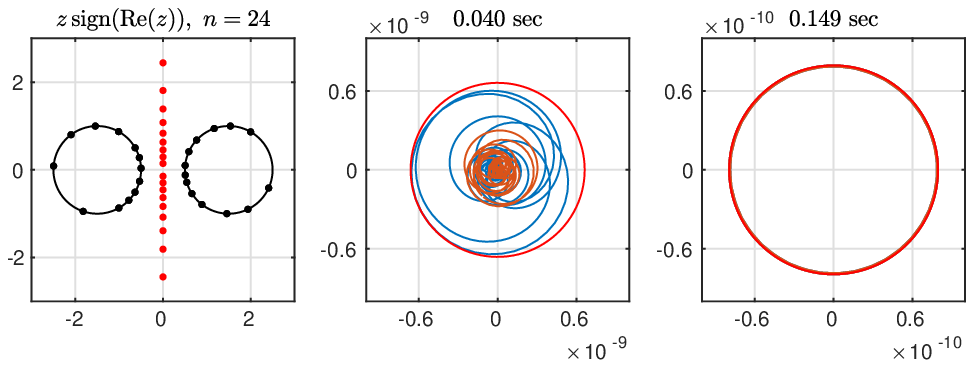}
\end{center}
\caption{\label{pair7}Approximations on the boundaries of
an annulus (doubly connected) and a union of disks (disconnected).
In each case the errors coresponding to the two disjoint
boundary components are plotted in different
colors.  For the annulus, the small red mark near the
origin reflects the fact that the 
best approximation has error $57.1$ times smaller on
the outer circle than the inner one.}
\end{figure}

Our final pair of complex examples, shown in Fig.~\ref{pair7}, involves
domains of more complicated connectivity.  The upper example approximates
the function $(1-z^{-2})^{1/2}$ on the boundary of the
annulus $1\le |z| \le 2$ (500 equispaced points on the outer circle together
with 500 points each in a {\tt tanh(12*linspace(-1,1))} distributions on the
upper and lower halves of the inner circle).  Note
that as usual, the poles cluster near the singularities on the boundary, which in
this case are at $\pm 1$.  The lower example approximates the
function $z\kern 1pt \hbox{sign}(\hbox{Re}(z))$ on the union of two circles
of radius $1$ about $-1.5$ and $1.5$ (1000 equispaced points
on each circle).  This function is not globally analytic, and
both AAA and AAA-Lawson tend to have difficulties with such problems.
Indeed if
$z\kern 1pt \hbox{sign}(\hbox{Re}(z))$ is replaced by
$\hbox{sign}(\hbox{Re}(z))$, the iteration fails.

Reviewing the 14 AAA-Lawson error curves displayed in
Figs.~\ref{expzpair}--\ref{pair7} (or error dots, in the case of
Fig.~\ref{randpair}), we note that it seems vividly apparent from
the near-maximal values at most of the points that a near-minimax
solution has been found.  Proving this would be a challenge,
however, although the bound $E^*\ge \min_{z\in Z} |r(z)-f(z)|$
follows from arguments related to Rouch\'e's theorem in cases
where the error curve is a near-circle of sufficiently high winding
number~\cite{gutk,htg,klotz,cf}.

\section{Numerical examples, real}
For real approximation on real intervals and unions of real
intervals, AAA-Lawson, like AAA itself, is less reliable
than in the complex case but retains its great speed and flexibility.
We shall present eight examples, grouping them again in pairs.

\begin{figure}
\begin{center}
\includegraphics{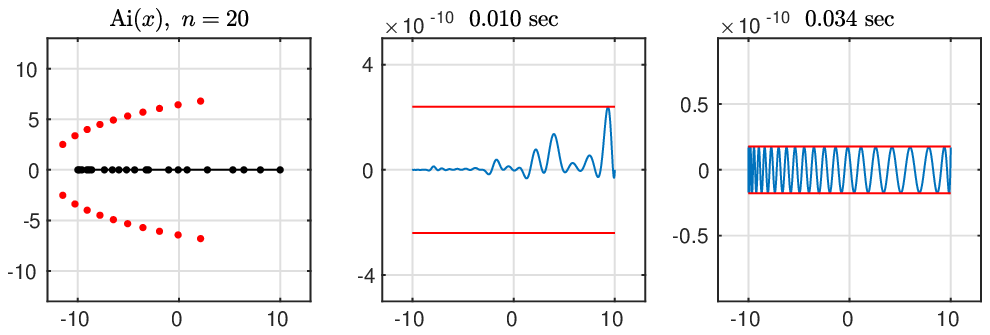}

\vskip .15in

\includegraphics{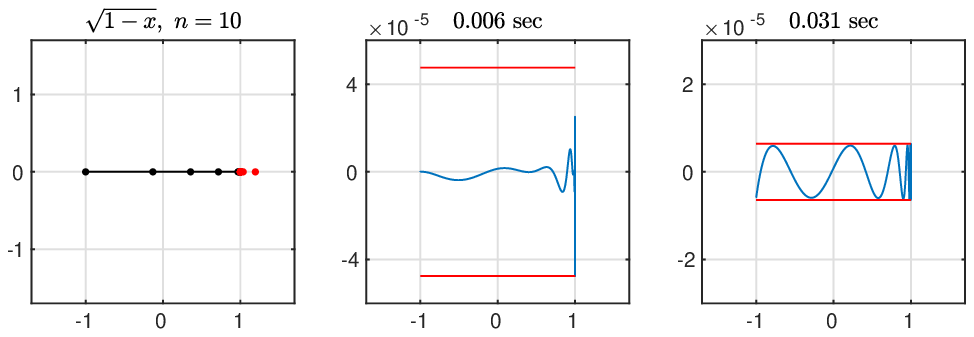}
\end{center}
\caption{\label{real1}Approximation on real intervals
of an analytic function and a 
function with a singularity at one endpoint.  The latter has exponentially clustered
poles approaching within a distance {\rm 1.4e-8} of $x=1$.}
\end{figure}

The first example of Fig.~\ref{real1} approximates $\hbox{Ai}(x)$
on $[-10,10]\kern .3pt $, which is discretized by 1000 points in a
Chebyshev distribution.  Note how the poles lie along curves in the
left half-plane, where the function is larger.  (The study of such
curves in approximation theory goes back to an investigation of
roots of Taylor polynomials of $e^z$ by Szeg\H o~\cite{szego}.)
In this case of an analytic function on a single interval,
Chebfun's {\tt minimax} gets the answer in 1.7 secs.\ and its
Carath\'eodory--Fej\'er command {\tt cf} does it in just 0.05
secs.~\cite{vandeun}.  The second example of the figure considers
$(1-x)^{1/2}$, which has a singularity at the right endpoint,
discretized on the grid {\tt tanh(linspace(-12,12,1000))}.\ \
The poles of this approximation cluster near $x=1$ at distances
15.3, 2.1, 0.19, 3.7e-2, 6.4e-3, 9.5e-4, 1.1e-4, 1.0e-5, 5.9e-7,
and 1.4e-8.  Chebfun {\tt minimax} is unsuccessful for this problem
with $n=10$, though it can handle degrees up to $n=7$.

\begin{figure}
\begin{center}
\includegraphics{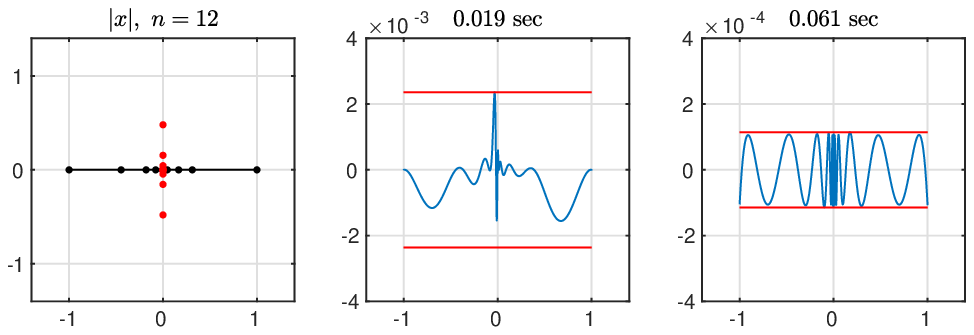}

\vskip .15in

\includegraphics{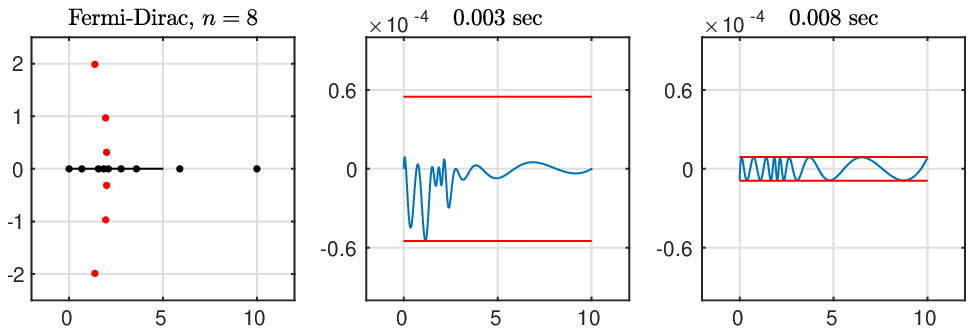}
\end{center}
\caption{\label{real2}Approximation of
functions with singularities (above) and near-singularities
(below) in the interval of approximation.}
\end{figure}

Figure~\ref{real2} turns to functions with a singularity or
near-singularity in the interior of the interval.  We pick
two examples where AAA and AAA-Lawson are successful, though
failures are common with problems of this kind.  The first example
is $|x|$, the problem made famous by Donald Newman, which is
discretized by transplants of {\tt tanh(linspace(-12,12))} to
both $[-1,0\kern .3pt ]$ and $[\kern .3pt 0,1]$; see~\cite{newman}
and~\cite[chapter 25]{atap}.  The AAA-Lawson error of 1.23e-4 is a
bit higher than the result 1.07e-4 computed by Chebfun {\tt minimax}
in 1.2 seconds.  As in Figs.~\ref{squarepair} and~\ref{scpair},
we see the 12 poles lining up along a branch cut; their locations
are approximately $\pm 0.00138\kern .3pt i$, $\pm 0.0102\kern
.3pt i$, $\pm 0.0448\kern .3pt i$, $\pm 0.155\kern .3pt i$, $\pm
0.4780\kern .3pt i,$ and $\pm 1.98\kern .3pt i$.  The second is the
Fermi--Dirac function $1/(1+\exp(10(x-2)))$ on the interval $[\kern
.3pt 0,10\kern .3pt ]$, as discussed for
example in~\cite{pexsi,moussa}, for which AAA-Lawson
gets an error of 9.09e-6.  Chebfun {\tt minimax} gets the better
value 8.77e-6 in 0.2 secs.\ and {\tt cf} does the same in $0.05$
secs.\ \ (A more robust computational strategy for Fermi--Dirac
functions is to first transplant $[\kern .3pt 0,\infty)$ to $[-1,1]$
by a M\"obius transformation~\cite{moussa,chebexp}.)

\begin{figure}
\begin{center}
\includegraphics{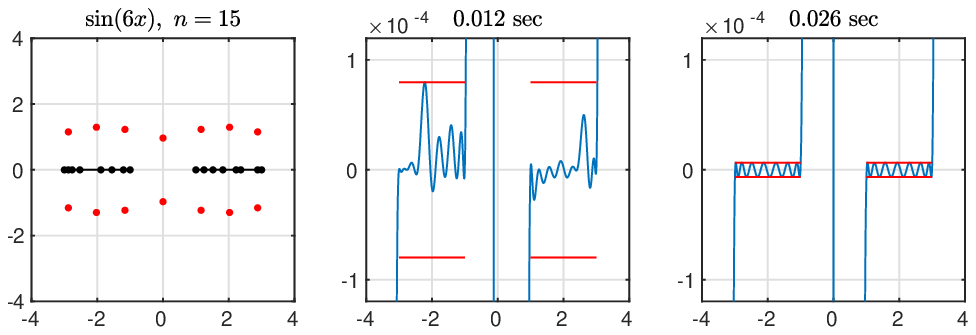}

\vskip .15in

\includegraphics{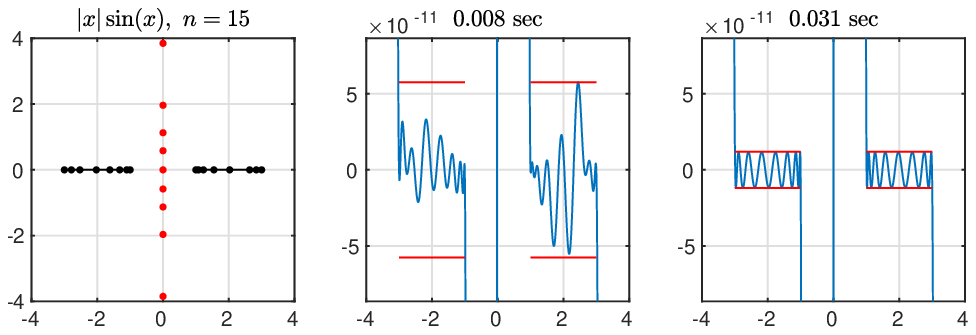}
\end{center}
\caption{\label{real3}Approximation of two functions each on a pair of disjoint
intervals.  The first example, $\sin(6x)$ is globally analytic, whereas the second, 
$|x|\sin(x)$, is analytic on each interval but not globally.  The very different
configurations of poles reflect this distinction.}
\end{figure}

Fig.~\ref{real3} considers a pair of problems on a union of two
intervals, $[-3,-1]\cup [1,3]$, each discretized by 500 points
in a Chebyshev distribution.  The first function, $\sin(6x)$, is
globally analytic, but the second, $|x|\sin(x)$, is not.
Note how the poles line up along the
imaginary axis, delineating once more an implicit branch cut.

\begin{figure}
\begin{center}
\includegraphics{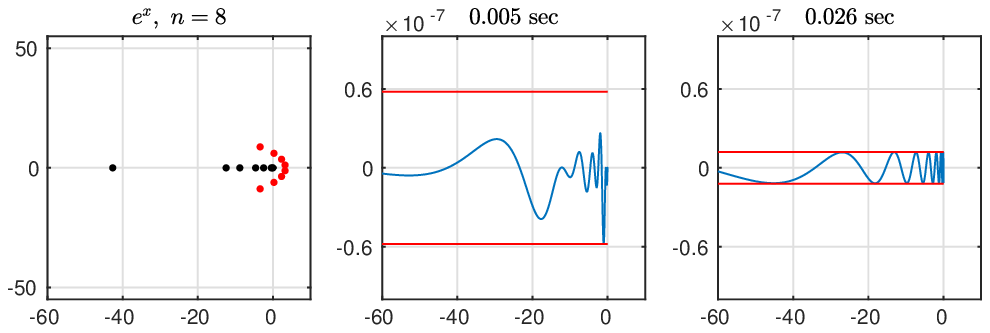}

\vskip .15in

\includegraphics{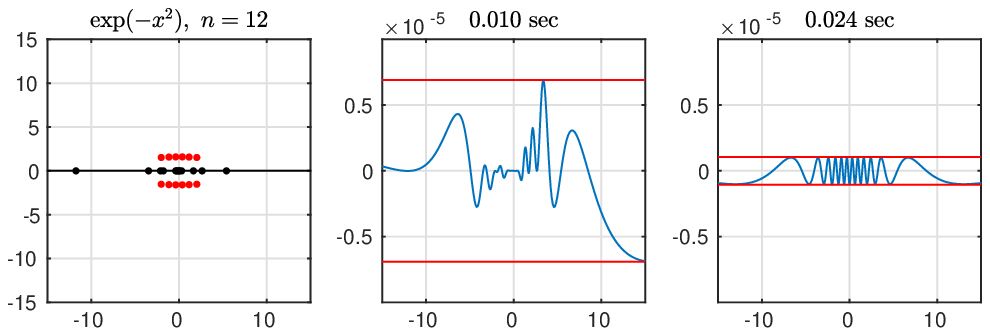}
\end{center}
\caption{\label{real4}Approximations
on $(-\infty,0\kern .3pt]$ and $(-\infty,\infty)$.}
\end{figure}

The final pair of examples, shown in Fig.~\ref{real4}, are posed
on infinite intervals.  The first is the Cody--Meinardus--Varga
problem of approximation of $e^x$ on $(-\infty,0\kern .3pt]$
\cite{cmv,aaa,talbot,atap}.  As described in section~4
of~\cite{talbot}, one can compute approximations here by
transplantation of $(-\infty,0\kern .3pt]$ to $[-1,1]$ followed
by CF or {\tt minimax} approximation, but here, we approximate
directly on the original untransplanted interval, which is
discretized by 2000 points logarithmically spaced from $-10^6$
to $-10^{-6}$.  The success of such a computation highlights
the extraordinary flexibility and stability of barycentric
representations based on support points selected by AAA.  The second
example of Fig.~\ref{real4} shows approximation of $\exp(-x^2)$ on
$(-\infty,\infty)$, discretized by 100 equispaced points in $[-1,1]$
concatenated with 500 logarithmically spaced points in $[1,10^6]$
and their negatives in $[-10^6,-1]$.  AAA-Lawson improves the error
from 6.92e-6 to 1.04e-6.

In general, we believe that the safest way to compute a real minimax
approximation on a real interval is usually by the Remez algorithm
as implemented in Chebfun {\tt minimax}, or if the function is
smooth, by CF approximation as implemented in {\tt cf}, in both cases
perhaps after softening up the problem by a M\"obius transformation.
The AAA-Lawson approach is most important in cases where these simple
tools are inapplicable, such as unbounded or disjoint intervals as
in Figs.~\ref{real3} or~\ref{real4}.

\section{\label{steps}Convergence properties}
Our experience with applying the AAA-Lawson method to hundreds of
examples can be summarized as follows: for analytic functions on
well resolved complex domains, it almost always converges, and for
nonanalytic functions or real domains, it often converges.  To give
more detail about our experience, here are summaries of the six
contexts we are aware of in which AAA-Lawson is most likely to fail.

{\em 1.~Discretization too coarse.}
Most applications involve discretization of a continuum, and trouble
often arises if the discretization is too coarse, especially near
singular points where poles need to accumulate.  Perhaps this is
unsurprising since even existence of best approximations fails in
general on discrete domains, as mentioned in section~\ref{theory}.
As indicated in the discussion of examples in the last two sections,
we routinely use Chebyshev-type sample point clustering near
nonsingular corners or endpoints of domains and more extreme {\tt
tanh(linspace(-12,12,npts))} type clustering near singular points.

{\em 2.~Too close to machine precision.}
By default, AAA delivers an approximation with accuracy close to
machine precision, and attempted Lawson iterations from such a point tend
to take on a random character, leading to failure.  Instead, in
standard double precision arithmetic, it is best to use AAA-Lawson
for approximations with errors down to $10^{-12}$ or $10^{-13}$
but not much smaller.  The Chebfun {\tt aaa} code reflects this
by running without Lawson if no degree is specified, e.g.\
\verb|aaa(F,Z)|, and with Lawson if a degree is specified, e.g.\
\verb|aaa(F,Z,'degree',10)|.  These defaults can be overridden
by specifying \verb|aaa(...,'lawson',nsteps)| in which case exactly
{\tt nsteps} Lawson steps are taken, and none at all if
$\hbox{\tt nsteps}=0$.  When we want accuracy to more digits
than delivered by default parameters, we
specify a large value of {\tt nsteps}.

{\em 3.~Degeneracy related to symmetry.}
Failure often occurs if one attempts a calculation that does not
respect the symmetry of the problem, where the mathematically correct
best approximation is degenerate.  For example, an attempt to compute
a degree $3$ best approximation to $\exp(z^2)$ on the unit disk
will fail, because the result should be of degree $2$. 
If the degree specification is changed to 2, the calculation succeeds.

{\em 4.~Lack of analyticity.}  The examples of the last two sections
illustrated that AAA-Lawson has little trouble with functions
meromorphic in a disk or an annulus.  Failures often occur in the
approximation of more deeply nonanalytic functions, however.  For
example, the example with 100 random points of Fig.~\ref{randpair}
fails if $f(z)$ is changed from $\tan(z)$ to $|z|$.

{\em 5.~Real domains.}  Failures are also common in
approximation of real functions on real domains.  As discussed
in~\cite{aaa}, such problems are difficult for AAA itself.

{\em 6.~Period-2 oscillations.}  Sometimes the Lawson iteration
enters into a cycle in which one pattern of weights and errors
appears at odd steps and another at even ones.
For example, this happens with the Fermi-Dirac example
of Fig.~\ref{real2} if $1/(1+\exp(10(x-2)))$ is changed to
$1/(1+\exp(50(x-2)))$, though the problem goes away if more sample
points are taken in the transition region.  In at least some
cases, convergence can also be recovered by
underrelaxation in the update formula (\ref{update}).

All of these failure modes reflect mathematical issues of substance
and point to interesting problems for future research.
As an engineering matter, the Chebfun code includes precautions to
minimize the risk of trouble arising from these sources, the most
basic of which is to revert to the AAA solution if AAA-Lawson fails
to make an improvement.  With time, we expect the engineering to
be further improved.

It is important to ask, what might be proved theoretically?
By making suitable assumptions, such as the use of exact arithmetic,
one could work around a number of the difficulties (1)--(6).
Still, proving anything
is far from straightforward, since even the theory of the linear
Lawson algorithm has encountered a number of obstacles, and here
we are working with a nonlinear barycentric variant.  We hope
that it may be possible to prove that for problems sufficiently
generic in an appropriate sense, and given a sufficiently close
initial guess and the use of a line search where necessary rather
than always taking the full update (\ref{update}), the iteration
is guaranteed to converge to a local minimum.  By analyzing the
perturbation properties of the SVD problem (\ref{barylin}), we have
made progress toward such a result, especially for real approximation
problems and complex problems with nearly-circular error curves.
However, we do not have a result comprehensive enough to report here.

\section{\label{sec-disc}Discussion}
The AAA-Lawson algorithm makes it easy for the first time
to compute
real and complex minimax rational approximations on all kinds of domains.
In Chebfun, for example, the commands
\begin{verbatim}

    Z = exp(2i*pi*(1:500)/500);
    F = exp(Z);
    r = aaa(F,Z,'degree',3);

\end{verbatim}
produce a function handle {\tt r} for the best degree
$3$ rational approximation of $f$ on the unit circle in a tenth
of a second on a laptop.  The calculation {\tt norm(F-r(Z),inf)}
then gives 9.9318e-6, matching the result published in~\cite{ew}
many years ago.

What makes the algorithm so effective is that it combines the
exceptional stability of barycentric rational representations, as
exploited by the AAA algorithm~\cite{aaa}, with the long-established technique
of IRLS iteration to improve the error to minimax---though in
a novel nonlinear barycentric context.  It is interesting that,
unlike its predecessors by Ellacott and Williams~\cite{ew} and
Istace and Thiran~\cite{istace}, AAA-Lawson is not based on an
attempt to satisfy optimality conditions.

As discussed in section~\ref{steps}, AAA-Lawson has little
theoretical foundation at present, and it also suffers from
just linear asymptotic convergence, sometimes at a low rate.
These drawbacks are not news!---as can be seen in this quote from
p.~50 of Osborne's book of 1985~\cite{osborne}:

\medskip

\begin{quotation}
\noindent The evidence presented here does not provide a recommendation 
for the technique [IRLS].\ \  It is shown that the
convergence rate is only first order in
general and that even this cannot be guaranteed.
\end{quotation}

\medskip

\noindent And yet, thanks to its simplicity and lack of dependence
on a characterization of optimal solutions, IRLS has enabled us to
develop an algorithm that is strikingly fast and robust.

Concerning the linear asymptotic convergence, it seems possible
that a method with improved behavior might be developed based
on Newton's method
combined with linear programming, tools applied effectively by Istace
and Thiran~\cite{istace} and Tang~\cite{tang}.  There are also
paradoxes to be investigated concerning the linear convergence
itself in problems with nearly-circular error curves, as in
Fig.~\ref{expzpair}.  Here, the nearly-constant error has the
effect that the Lawson weight distribution virtually stops changing from step to
step~\cite{ugrad}, and in particular, effectively never get close
to the sum of delta functions form that an asymptotic analysis is likely
to look for.  Nevertheless, the approximations in such cases often
converge quickly, and to add to the mystery, they converge much
faster still if (\ref{update}) is modified to depend on $|e_j|^2$
instead of $|e_j|$, although in other cases this
modifcation results in failure.  Perhaps an understanding of such effects might
lead to improvements in the AAA-Lawson algorithm even in cases with
error curves that are not nearly circular.

This article has considered only standard minimax approximations, without
weight functions.  Nonconstant weights are easily introduced by modifying
$(\ref{update})$.
Another restriction is that we have
treated only rational approximations of type
$(n,n)$, not type $(m,n)$ with $m\ne n$.  The more general problem
is certainly interesting, and AAA itself can be generalized to
$m\ne n$ as described in~\cite{aaa}.  However, though the ``Walsh
table'' of approximations of a function of all rational types is
fascinating, the overwhelming majority of applications are concerned
with types $(n,n)$ or $(n-1,n)$.

\section*{Acknowledgments}
We are grateful to Silviu Filip and Abi Gopal for assistance with
both software and mathematics.

\end{document}